\documentstyle[leqno,11pt]{article}

\newtheorem{theorem}{Theorem}[section]

\def\N{{\bf N}}

\def\1{{\bf 1}}

\begin{document}

\title{On the asymptotic densities of certain subsets of $\N^k$}
\author{{\sc L\'aszl\'o T\'oth}\\
University of P\'ecs,\\
Institute of Mathematics and Informatics,\\
Hungary, 7624 P\'ecs, Ifj\'us\'ag u. 6,\\
E-mail: {\tt ltoth@math.ttk.pte.hu}}
\date{\sl{Riv. Mat. Univ. Parma (6)} {\bf 4} (2001), 121-131}
\maketitle

\begin{abstract}
We determine the asymptotic density $\delta_k$
of the set of ordered $k$-tuples $(n_1,...,n_k)\in \N^k, k\ge 2$, such
that there exists no prime power $p^a$, $a\ge 1$, appearing in the canonical
factorization of each $n_i$, $1\le i\le k$, and
deduce asymptotic formulae with error terms
regarding this problem and analogous ones. We give numerical
approximations of the constants $\delta_k$ and improve the error term
of formula (1.2) due to {\sc E. Cohen}.

We point out that our treatment, based on certain
inversion functions, is applicable also in case $k=1$ in order to
establish asymptotic formulae with error terms regarding the densities
of subsets of $\N$ with
additional multiplicative properties. These generalize an often cited result
of {\sc G. J. Rieger}.

\end{abstract}



\section{Introduction}

Let $k\ge 2$ be a fixed integer. What is the asymptotic density $\delta_k$
of the set of ordered $k$-tuples $(n_1,...,n_k)\in \N^k$, such
that there exists no prime power $p^a$, $a\ge 1$, appearing in the canonical
factorization of each $n_i$, $1\le i\le k$ ?

This problem is analogous to the following one:
What is the asymptotic density $d_k$
of the set of $k$-tuples which are relatively prime, i.e. $k$-tuples
$(n_1,...,n_k)\in \N^k$ such that
there exists no prime $p$, appearing in the canonical
factorization of each $n_i$, $1\le i\le k$ ?

It is known that $d_k= 1/\zeta(k)$, where $\zeta$ is the Riemann
zeta function, and this value can be considered as
the probability that $k$ integers ($k\ge 2$) chosen
at random are relatively prime.  More precisely,
$$
N_k(x): = \# \{(n_1,...,n_k)\in (\N \cap [1,x])^k : gcd
(n_1,...,n_k)=1\} = \frac1{\zeta(k)}x^k + R_k(x), \leqno (1.1)
$$
where $R_k(x)=O(x^{k-1})$ for $k>2$, $R_2(x)=O(x\log x)$ for $k=2$, and
$d_k=\lim_{x\to \infty} N_k(x)/x^k$ $=1/\zeta(k)$.
This result goes back to the work of {\sc J. J. Sylvester} \cite{S}
and {\sc D. N. Lehmer} \cite{L}, see also {\sc J. E. Nymann} \cite{N}.

There are several generalizations of (1.1) in the literature. For example,
let $S$ be an arbitrary subset of $N$. Then
$$
N_k(x,S): = \# \{(n_1,...,n_k) \in (\N \cap [1,x])^k : gcd
(n_1,...,n_k)\in S\} = \frac{\zeta_S(k)}{\zeta(k)}x^k + T_k(x), \leqno (1.2)
$$
where $$\zeta_S(k)=\sum_{n=1 \atop{n\in S}}^{\infty} \frac1{n^k}$$
and $T_k(x)=O(x^{k-1})$ for $k>2$, $T_2(x)=O(x\log^2 x)$ for $k=2$,
for every $S\subseteq \N$, due to {\sc E. Cohen} \cite{C}.
Therefore the asymptotic density
of the set of ordered $k$-tuples $(n_1,...,n_k)$ for which
$gcd (n_1,...,n_k)$ belongs to $S$ is
$\lim_{x\to \infty} N_k(x,S)/x^k =\frac{\zeta_S(k)}{\zeta(k)}$.

{\sc J. E. Nymann} \cite{Nii} shows that if the characteristic function
$\rho_S$ of $\emptyset \ne S\subseteq \N$ is completely multiplicative and if
$\# \{n: n\in S\cap [1,x]\}= A x + O(1)$, where $A$ is the asymptotic
density of $S$, then
$$
\# \{(n_1,...,n_k)\in (S \cap [1,x])^k : gcd
(n_1,...,n_k)=1\} = A^k \prod_{p\in S} \left(1-\frac1{p^k}\right) x^k +
R_k(x),
\leqno (1.3)
$$
where $R_k(x)$ is the same as above.
Therefore, if $P_k^S(n)$ denotes the probability that $k$ integers ($k\ge 2$)
chosen
at random from $S\cap [1,n]$ are relatively prime, then
$$
\lim_{n\to \infty} P_k^S(n)= \prod_{p\in S} \left(1-\frac1{p^k}\right).
$$
This result can be applied for $S=\{ n: \textstyle{ gcd } (n, p_1\cdots p_r)=1
\}$, where $\{p_1,...,p_r\}$ is a given finite set of distinct primes.

Now return to the problem at the beginning. It is obvious that $\delta_k \ge d_k=
1/\zeta(k)$ for every $k\ge 2$ and thus $\lim_{k\to \infty} \delta_k=1$. Which
is the exact value of $\delta_k$ ?

In order to solve this problem we use the concept of the unitary divisor.
For $d,n \in \N$, $d$ is called a unitary divisor (or block divisor) of $n$
if $d|n$ and $gcd (d,n/d)=1$, notation $d||n$.
Various other problems concerning unitary divisors, including properties of
arithmetical functions and arithmetical convolutions defined by unitary divisors,
have been studied extensively in the literature, see for example \cite{M} and
its bibliography. Denote the greatest common unitary divisor
of $n_1,...,n_k$ by $gcud (n_1,...,n_k)$.

Our question can be reformulated in this way:
What is the asymptotic density $\delta_k$
of the set of ordered $k$-tuples $(n_1,...,n_k)$ such that
$gcud (n_1,...,n_k)=1$, or more generally, $gcud (n_1,...,n_k)\in S$ ?

Furthermore, what is the probability that for $k$ integers $n_1,...,n_k$
chosen at random from $S\cap [1,n]$ one has $gcud (n_1,...,n_k)=1$ ?

In this paper we determine the value $\delta_k$ and deduce asymptotic formulae
with error termes analogous to (1.1) -(1.3), regarding these problems.
We give numerical
approximations of the constants $\delta_k$ and also improve the error term
of (1.2) of {\sc E. Cohen}.

The treatment we use is based on on the inversion functions $\mu^*_S$ and
$\mu_S$ attached to the subset $S$.
We point out that this is applicable also in case $k=1$ in order to
establish asymptotics regarding the densities of certain subsets $S$ of $\N$,
generalizing in this way an often cited result of {\sc G. J. Rieger} \cite{R}.

Note that the value $\delta_2$ is given by
{\sc D. Suryanarayana} and {\sc M. V. Subbarao} \cite{SS}, Corollary 3.6.3,
applying other arguments as those of the present paper.

Using the concept of regular cross-convolution, see \cite{Ti}, \cite{Tviii},
it is possible to deduce more general results, including (1.1) - (1.3)
and (2.1) and (2.4) of this paper. We do not go into details.

\section{Results}

Let $S\subseteq \N$. We say that $S$ is (completely) multiplicative if
$1\in S$ and its characteristic function $\rho_S(n)$ is (completely)
multiplicative.  Define the function $\mu^*_S(n)$ by
$$
\sum_{d||n} \mu^*_S(d)=\rho_S(n), \quad n\in \N,
$$
that is
$$
\mu^*_S(n)= \sum_{d||n} \rho_S(d)\mu^*(n/d), \quad n\in \N,
$$
where the sums are extended over the unitary divisors of $n$ and
$\mu^*(n): = \mu^*_{\{1\}}(n)=(-1)^{\omega(n)}$, $\omega(n)$ denoting
the number of distinct prime factors of $n$.

Furthermore, let $\phi(n)$ and $\theta(n)$ denote Euler's
function and the number of squarefree divisors of $n$, respectively.

\begin{theorem}  If $k\ge 2$ and $S$ is an arbitrary subset of $\N$, then
$$
\# \{(n_1,...,n_k)\in (\N \cap [1,x])^k : gcud
(n_1,...,n_k)\in S\} = \delta_k(S) x^k + V_k(x,S), \leqno (2.1)
$$
where
$$
\delta_k(S)=\sum_{n=1}^{\infty} \frac{\mu^*_S(n)\phi^k(n)}{n^{2k}}
$$
and the remainder term can be evaluated as follows:

{\rm (1)} $V_k(x,S)=O(x^{k-1})$ for $k>2$ and for an arbitrary $S$,

{\rm (2)} $V_2(x,S)=O(x\log^4 x)$ for an arbitrary $S$,

{\rm (3)} $V_2(x,S)=O(x\log^2 x)$ for an $S$ such that
$\sum_{n\in S} \frac{\theta(n)}{n}< \infty$ (in particular for every finite $S$)
and for every multiplicative $S$,

{\rm (4)} $V_2(x,S)=O(x)$ for every multiplicative $S$ such that
$\sum_{p\notin S} \frac1{p}< \infty$ (in particular if the set
$\{p: p\notin S\}$ is finite).

If $S$ is multiplicative, then
$$
\delta_k(S)= \prod_p \left(1-(1-\frac1{p})^k \sum_{a=1\atop{p^a\notin S}}^{\infty}
\frac1{p^{ak}}\right).
$$
If $S=\{1\}$, then
$$
\delta_k:= \delta_k(\{1\}) = \prod_p \left(1- \frac{(p-1)^k}{p^k(p^k-1)}\right).
$$
\end{theorem}
\vskip2mm

\begin{theorem}  If $k\ge 2$ and $S$ is an arbitrary subset of $\N$, then
the asymptotic densities of the sets of ordered $k$-tuples $(n_1,...,n_k)$ such
that $gcud (n_1,...,n_k)\in S$ and $gcud (n_1,...,n_k)=1$ are
$\delta_k(S)$ and $\delta_k$, respectively, given in {\rm Theorem 2.1}.
\end{theorem}

\begin{theorem} Let $p_n$ denote the $n$-th prime and for $r\in \N$ let
$N=10^r/2$. Then
$$
\delta_k \approx \prod_{n=1}^N \left(1- \frac{(p_n-1)^k}{p_n^k(p_n^k-1)}\right)
$$
is an approximation of $\delta_k$ with $r$ exact decimals.

In particular, $\delta_2 \approx 0.8073,  \delta_3 \approx 0.9637,
\delta_4 \approx 0.9924, \delta_5 \approx 0.9983,
\delta_6 \approx 0.9996, \delta_7 \approx 0.9999$, with
$r=4$ exact decimals.
\end{theorem}

\begin{theorem} For $k=2$ the error term $R_2(x)$ of {\rm (1.2)}
can be improved into $R(x,S)$, where

{\rm (i)} $R(x,S)=O(x\log x)$ for an $S$ such that
$\sum_{n\in S} \frac1{n}< \infty$ (in particular for every finite $S$)
and for every multiplicative $S$,

{\rm (ii)} $R(x,S)=O(x)$ for every multiplicative $S$ such that
$\sum_{p\notin S} \frac1{p}< \infty$ (in particular if the set
$\{p: p\notin S\}$ is finite).
\end{theorem}

{\bf Remark.} It is noted in [1] that if $k=2$ and if the function $\mu_S(n)$
is bounded, cf. proof of Theorem 2.4 of the present paper, then the error
term is $R_2(x)=O(x \log x)$.

\begin{theorem} Suppose that $S\subseteq \N$ is multiplicative and
$\min \{a: p^a\notin S\}\ge r \ge 2$ for every prime $p$. Then
$$
\sum_{n\le x} \rho_S(n) =d(S)x +O(\sqrt[r]{x}), \leqno (2.2)
$$
where
$$
d(S)=\prod_{p} \left(1-\frac1{p}\right) \left( 1+ \sum_{a=1
\atop{p^a\in S}}^{\infty} \frac1{p^a}\right)   \leqno (2.3)
$$
is the asymptotic density of $S$.
\end{theorem}

{\bf Remark.} In the special case $S=$ the set of $K$-void integers we reobtain
from (2.2) the result of {\sc G. J. Rieger} \cite{R}.
The $K$-void integers are defined as follows. Let
$K$ be a nonempty subset of $\N\setminus \{1\}$. The number $n$ is called
$K$-void if $n=1$ or $n>1$ and
there is no prime power $p^a$, with $a\in K$,
appearing in the canonical factorization of $n$.

Does the density exist for an arbitrary multiplicative subset $S$ ?
Yes, and it is $d(S)$ given by (2.3), where the infinite product is considered
to be $0$
when it diverges (if and only if $\sum_{p\notin S} \frac1{p}= \infty$).
This follows from a well-known result of {\sc E. Wirsing} \cite{W} concerning
the mean-values of certain multiplicative functions $f$.
A short direct proof for the case $f$ multiplicative and $0\le f(n)\le 1$
for $n\ge 1$, hence applicable for the characteristic function of an arbitarary
multiplicative $S$, is given in the book of {\sc G. Tenenbaum}, \cite{T}, p. 48.

\begin{theorem} Let $k\ge 2$ and suppose that $S$ is a completely
multiplicative subset of $\N$
such that $\# \{n: n\in S\cap [1,x]\}= A x + O(1)$. Then
$$
\# \{(n_1,...,n_k)\in (S \cap [1,x])^k : gcud
(n_1,...,n_k)=1\} = A^k \beta_k(S) x^k + T_k(x),
\leqno (2.4)
$$
where
$$
\beta_k(S)= \prod_{p\in S} \left(1- \frac{(p-1)^k}{p^k(p^k-1)}\right),
$$
and $T_k(x)=O(x^{k-1})$ for $k>2$, $T_2(x)=O(x\log^2 x)$ for $k=2$.

If $Q_k^S(n)$ denotes the probability that for $k$ integers $n_1,...,n_k$
chosen at random from $S\cap [1,n]$ one has $gcud (n_1,...,n_k)=1$, then
$$
\lim_{n\to \infty} Q_k^S(n)= \beta_k(S).
$$
\end{theorem}

\section{Proofs}

{\bf Proof of Theorem 2.1}
Using the definition of $\mu_S^*$,
the fact that $d|| gcud (n_1,...,n_k)$ if and only if $d||n_i$ for every
$1\le i \le k$, which can be checked easily, and the well-known estimate
$$
\sum_{n\le x \atop{gcd (n,m)=1}} 1 = \frac{\phi(m)x}{m}
+ O(\theta(m))
$$
which holds uniformly for $x\ge 1$ and $m\in \N$, we obtain
$$
\# \{(n_1,...,n_k)\in (\N \cap [1,x])^k : gcud (n_1,...,n_k)\in S\} =
 \sum_{n_1,...,n_k\le x} \rho_S(gcud (n_1,...,n_k))=
$$
$$
= \sum_{n_1,...,n_k\le x} \sum_{d||(n_1,...,n_k)}
\mu^*_S(d) =
\sum_{n_1,...,n_k\le x} \sum_{d||n_1,...,d||n_k} \mu^*_S(d)=
$$
$$
=\sum_{d\le x} \mu^*_S(d) \sum_{a_i\le x/d \atop{(a_i,d)=1
\atop{1\le i\le k}}} 1
=\sum_{d\le x} \mu^*_S(d) \left(\sum_{a\le x/d \atop{(a,d)=1}} 1 \right)^k =
$$
$$
=\sum_{d\le x} \mu^*_S(d) \left(\frac{x\phi(d)}{d^2}+O(\theta(d))\right)^k
=\sum_{d\le x} \mu^*_S(d) \left(\frac{x^k\phi^k(d)}{d^{2k}} +
O(\frac{x^{k-1}\theta(d)}{d^{k-1}})\right)=
$$
$$
=x^k \sum_{d\le x} \frac{\mu^*_S(d)\phi^k(d)}{d^{2k}} +
O\left(x^{k-1} \sum_{d\le x} \frac{|\mu^*_S(d)|\theta(d)}{d^{k-1}}\right)=
$$ $$
=\delta_k(S) x^k + O\left( x^k \sum_{d>x} \frac{|\mu^*_S(d)|}{d^k}\right) +
O\left(x^{k-1} \sum_{d\le x} \frac{|\mu^*_S(d)|\theta(d)}{d^{k-1}}\right).
$$

The given error term yields now from the next statements:

(a)
$$ \sum_{n\le x} \frac{\theta(n)}{n^s}=\cases
{O(\log^2 x), & $s=1$, \cr  O(1), & $s>1$. \cr}
$$ $$
\sum_{n\le x} \frac{\theta^2(n)}{n^s}=\cases
{O(\log^4 x), & $s=1$, \cr O(1), & $s>1$, \cr}
$$ $$ \sum_{n>x} \frac1{n^s}= O(\frac1{x^{s-1}}), \qquad
\sum_{n>x} \frac{\theta(n)}{n^s}=O(\frac{\log x}{x^{s-1}}), \quad s>1.
$$

(b) For an arbitrary $S\subseteq \N$ and for every $n\in \N$,
$|\mu^*_S(n)|\le \sum_{d||n} \rho_S(d)
\le \theta(n)$, $|\mu^*_S(n)|\theta(n)\le \sum_{d||n} \rho_S(d)
\theta(d)\theta(n/d)$ and
$$
\sum_{n\le x} \frac{|\mu^*_S(n)|\theta(n)}{n}\le \sum_{d\le x}
\frac{\rho_S(d)\theta(d)}{d}
\sum_{e\le x/d} \frac{\theta(e)}{e}= $$ $$= O\left(\log^2 x \sum_{d\le x}
\frac{\rho_S(d)\theta(d)}{d}\right)
= \cases{O(\log^2 x), &if $\sum_{n=1}^{\infty} \frac{\rho_S(n)\theta(n)}{n}
           < \infty$,\cr
         O(\log^4 x), &otherwise. \cr}
$$

(c) If $S$ is multiplicative, then $\mu^*_S$ is multiplicative too,
$\mu^*_S(p^a)=\rho_S(p^a)-1$ for every prime power $p^a$ ($a\ge 1$)
and $\mu^*_S(n)\in \{-1,0,1\}$ for each $n\in \N$.

(d) Suppose $S$ is multiplicative. Then
$$
\sum_p \sum_{a=1}^{\infty} \frac{|\mu^*_S(p^a)|\theta(p^a)}{p^a} =
2 \sum_p \sum_{a=1}^{\infty} \frac{1-\rho_S(p^a)}{p^a} \le
$$ $$
\le 2\sum_p \left(\frac{1-\rho_S(p)}{p}+\sum_{a=2}^{\infty} \frac1{p^a}\right)
= 2 \sum_{p\in S} \frac1{p(p-1)} + 2 \sum_{p\notin S} \frac1{p-1}\le
$$
$$
\le 4\left(\sum_{p\in S} \frac1{p^2} + \sum_{p\notin S} \frac1{p}\right) <
\infty \quad
{\rm if} \quad \sum_{p\notin S} \frac1{p} < \infty.
$$
It follows that in this case the series $\sum_{n=1}^{\infty}
\frac{|\mu^*_S(n)|\theta(n)}{n}$ is convergent.

If $S$ is multiplicative, then the series giving  $\delta_k(S)$ can be
 expanded into an infinite product of Euler-type.

{\bf Proof of Theorem 2.2} This is a direct consequence of Theorem 2.1.

{\bf Proof of Theorem 2.3}
Consider the series of positive terms
$$
\sum_p \log \left(1-\frac{(p-1)^k}{p^k(p^k-1)}\right)^{-1}=
$$
$$
=\sum_{n=1}^{\infty} \log \left(1+\frac{(p_n-1)^k}{p_n^k(p_n^k-1)-(p_n-1)^k}
\right) =  - \log \delta_k,
$$
where $p_n$ denotes the $n$-th prime.

The $N$-th order error $R_N$ of this series can be evaluated
as follows:
$$
R_N: = \sum_{n=N+1}^{\infty} \log \left(1+\frac{(p_n-1)^k}{p_n^k(p_n^k-1)-
(p_n-1)^k} \right)
< \sum_{n=N+1}^{\infty} \frac{(p_n-1)^k}{p_n^k(p_n^k-1)-(p_n-1)^k}<
$$ $$
< \sum_{n=N+1}^{\infty} \frac1{p_n^k-1}\le
 \sum_{n=N+1}^{\infty} \frac1{p_n^2-1}.
$$

Now using that $p_n>2n$, valid for $n\ge 5$, we have
$$
R_N< \sum_{n=N+1}^{\infty} \frac1{4n^2-1}=\frac1{2}
\sum_{n=N+1}^{\infty} \left( \frac1{2n-1} - \frac1{2n+1}\right) =\frac1{2(2N+1)}.
$$
In order to obtain an approximation with $r$ exact decimals we use the
condition
$$
\frac1{2(2N+1)}\le \frac1{2}\cdot 10^{-r}
$$
and obtain $N\ge \frac1{2} ( 10^r -1)$.

The numerical values were obtained using the software package MAPLE.

{\bf Proof of Theorem 2.4}
Define the function $\mu_S(n)$ by
$$
\sum_{d|n} \mu_S(d)=\rho_S(n), \quad n\in \N,
$$
that is
$$
\mu_S(n)= \sum_{d|n} \rho_S(d)\mu(n/d), \quad n\in \N,
$$
where $\mu(n): = \mu_{\{1\}}(n)$ is the M\"obius function, see \cite{C}.
We have
$$
N_k(x,S): =\# \{(n_1,...,n_k)\in (\N \cap [1,x])^k : gcd (n_1,...,n_k)\in S\} =
$$
$$
= \sum_{n_1,...,n_k\le x} \rho_S(gcd (n_1,...,n_k))=
\sum_{n_1,...,n_k\le x} \sum_{d|(n_1,...,n_k)} \mu_S(d)
$$
and the proof runs parallel to the proof of Theorem 2.1.

{\bf Proof of Theorem 2.5}
$$
N_1(x,S) = \sum_{n\le x} \rho_S(n) = \sum_{n\le x} \sum_{d|n} \mu_S(d)
=x \sum_{d\le x} \frac{\mu_S(d)}{d} + O( \sum_{d\le x} |\mu_S(d)|).
$$
Here $\mu_S$ is multiplicative, $\mu_S(p^a)=\rho_S(p^a)-\rho_S(p^{a-1}),
a\ge 1$ and since $p,p^2,...,p^{r-1}\in S$ we have $\mu_S(p)=\mu_S(p^2)=
... =\mu_S(p^{r-1})=0$ for every prime $p$. Hence for each $n\in \N$,
$|\mu_S(n)|\le \rho_{L_r}(n)$,
where $L_r$ is the set of $r$-full numbers, i. e.
$L_r = \{1\} \cup \{n>1: p|n \Rightarrow p^r| n \}$.
We get
$$
N_1(x,S) = d(S) x + O(x \sum_{d>x} \frac{\rho_{L_r}(d)}{d}) +
O(\sum_{d\le x}\rho_{L_r}(d)),
$$
and using the elementary estimate
$$
\sum_{n\le x} \rho_{L_r}(n)= C \sqrt[r]{x} + O(\sqrt[r+1]{x}),
$$
where $C$ is a positive constant, due to {\sc P. Erd\H os} and
{\sc G. Szekeres} \cite{ESz}, obtain the given result.

{\bf Proof of Theorem 2.6}
$$
\# \{(n_1,...,n_k)\in (S \cap [1,x])^k : gcud (n_1,...,n_k)=1\}=
$$
$$
= \sum_{n_1\le x} \rho_S(n_1)... \sum_{n_k\le x} \rho_S(n_k)
\sum_{d||(n_1,...,n_k)} \mu^*(d) =
\sum_{n_1\le x} \rho_S(n_1)... \sum_{n_k\le x} \rho_S(n_k)
\sum_{d||n_1,...,d||n_k} \mu^*(d)=
$$
$$
=\sum_{d\le x} \mu^*(d) \sum_{a_1\le x/d \atop{(a_1,d)=1}}\rho_S(da_1)
... \sum_{a_k\le x/d \atop{(a_k,d)=1}} \rho_S(da_k)=
$$  $$
= \sum_{d\le x} \rho_S(d) \mu^*(d) \left(\sum_{a\le x/d \atop{(a,d)=1}}
\rho_S(a)\right)^k.
$$

Here we use the estimate, valid for every $\ell \in \N$,
$$
\sum_{n\le x \atop{gcd (n,\ell)=1}} \rho_S(n)=
\sum_{n\le x} \rho_S(n) \sum_{d| gcd (n,\ell)} \mu(d)=
$$
$$
=\sum_{d|\ell} \mu(d)\rho_S(d) \sum_{e\le x/d} \rho_S(e)
=\sum_{d|\ell} \mu(d)\rho_S(d) \left(A\frac{x}{d}+O(1) \right) =
$$
$$
= A x \prod_{p|\ell \atop{p\in S}} (1-\frac1{p}) +O(\theta(\ell))
$$
and obtain the desired result, see the proof of Theorem 2.1.


\begin{thebibliography}{99}

\bibitem{C} {\sc E. Cohen}, {\it Arithmetical functions associated with
arbitrary sets of integers}, Acta Arith., {\bf 5} (1959), 407-415.

\bibitem{ESz} {\sc P. Erd\H os, G. Szekeres}, {\it \"Uber die Anzahl der
Abelschen Gruppen gegebener Ordnung und \"uber ein verwandtes
zahlentheoretisches Problem}, Acta Sci. Math. (Szeged), {\bf 7} (1935), 95-102.

\bibitem{L} {\sc D. N. Lehmer}, {\it An asymptotic evaluation of certain totient
sums}, Amer. J. Math., {\bf 22} (1900), 293-355.

\bibitem{M} {\sc P. J. McCarthy}, {\it Introduction to Arithmetical
Functions}, Springer Verlag, New York - Berlin - Heidelberg - Tokyo, 1986.

\bibitem{N} {\sc J. E. Nymann}, {\it On the probability that $k$ positive
integers are relatively prime}, J. Number Theory, {\bf 4} (1972), 469-473.

\bibitem{Nii} {\sc J. E. Nymann}, {\it On the probability that $k$ positive
integers are relatively prime, II.}, J. Number Theory, {\bf 7} (1975), 406-412.

\bibitem{R} {\sc G. J. Rieger}, {\it Einige Verteilungsfragen mit $K$-leeren
Zahlen, $r$-Zahlen und Primzahlen}, J. Reine Angew. Math., {\bf 262/263} (1973),
189-193.

\bibitem{SS} {\sc D. Suryanarayana, M. V. Subbarao}, {\it Arithmetical functions
associated with the bi-unitary $k$-ary divisors of an integer},
Indian J. Math., {\bf 22} (1980), 281-298.

\bibitem{S} {\sc J. J. Sylvester}, {\it The Collected Mathematical Papers of
James Joseph Sylvester, vol. III.}, Cambridge Univ. Press, London - New York,
1909.

\bibitem{T} {\sc G. Tenenbaum}, Introduction to Analytic and Probabilistic
Number Theory, Cambridge  Univ. Press, 1995.

\bibitem{Ti} {\sc L. T\'oth}, {\it Asymptotic formulae concerning arithmetical
functions defined by cross-convolutions, I. Divisor-sum functions and
Euler-type functions}, Publ. Math. Debrecen, {\bf 50} (1997), 159-176.

\bibitem{Tviii} {\sc L. T\'oth}, {\it Asymptotic formulae concerning
arithmetical functions defined by
cross-convolutions, VIII. On the product and the quotient of $\sigma_{A,s}$
and $\phi_{A,s}$}, Riv. Mat. Univ. Parma (6), {\bf 2} (1999), 199-206.

\bibitem{W} {\sc E. Wirsing}, {\it Das aymptotische Verhalten von Summen \"uber
multiplikative Funktionen, II.}, Acta Math. Acad. Sci. Hung., {\bf 18}
(1967), 411-467.

\end{thebibliography}
\end{document}